\documentclass{math}
\usepackage{hyperref}
\hypersetup{
colorlinks=true,
urlcolor=blue,
citecolor=blue}
\usepackage[all]{xy}
\usepackage{amsfonts,amssymb,amsmath,amsgen,amsopn,amsbsy,theorem,graphicx,epsfig}
\usepackage{eufrak,amscd,bezier,latexsym,mathrsfs,enumerate,multirow}
\usepackage[utf8]{inputenc}\usepackage[english]{babel}
\usepackage[dvipsnames]{xcolor}
\usepackage[pagewise]{lineno}
\linenumbers

\yil{}
\vol{}
\fpage{}
\lpage{}
\doi{}

\title{Group invariant solutions and Conservation laws of the nonlinear Gardner-Kawahara equation. \\
}

\author[Sradharam Swain and Bikash Sahoo and Manjit Singh]{
\textbf{Sradharam Swain$^{1}$\thanks{Correspondence: sradharam1000@gmail.com}, Bikash Sahoo$^{2}$ Manjit Singh$^{3}$}\\
\\ 
$^{1}$Department of Mathematics, National Institute of Technology, Rourkela, India, \\
ORCID iD:https://orcid.org/0000-0001-5963-7856, contact-sradharam1000@gmail.com\\
\\
$^{2}$Department of Mathematics, National Institute of Technology, Rourkela, India, \\
 ORCID iD: https://orcid.org/0000-0003-0097-0915, contact-bikashsahoo@nitrkl.ac.in\\
 \\
 $^{3}$Department of Mathematics, Yadavindra College of Engineering Punjabi University \\of Campus Talwandi Sabo. contact-manjitcsir@gmail.com\\ ORCID iD: https://orcid.org/0000-0002-3988-4398
 
\\
\textcolor{red}{Sradharam Swain.}
\\ [1.8em]

\rec{.201}
\acc{.201}
\finv{..201}
}

\amssayisi{2010 {\itshape AMS Mathematics Subject Classification:} 76M60, 35C05, 70H33
\newline
\vspace{-15mm}
\begin{center}
		 \raisebox{-17ex}[0ex][0ex]{~~ \raisebox{.5ex}[0ex][0ex]{\footnotesize  This work is licensed under a Creative Commons Attribution 4.0 International License.}}
 		    \end{center}}

\newcommand{\bc}{\begin{center}}
\newcommand{\ec}{\end{center}}

\numberwithin{equation}{section}

\renewcommand{\phi}{\varphi}

\begin{document}

\maketitle

\begin{abstract}
	The present article studies the potential form of the nonlinear Gardner-Kawahara equation through the perspective of Lie symmetry analysis. Lie symmetry analysis was used to investigate abundant group-invariant solutions of the nonlinear Gardner-Kawahara equation.
	This method is used to provide geometric vector fields, as well as their commutative and adjoint relations. In this article, we have obtained the exact solution of the nonlinear Gardner-Kawahara equation in explicit form by different significant methods. Numerical simulation is used to explain the physical relevance of invariant solutions in 3D and 2D graphs.. Finally, by the conservation law multiplier, the complete set of local conservation laws of the equation for the arbitrary constant coefficients is well constructed with a detailed derivation. The conserved currents discovered in this study can help us better comprehend some of the physical processes that the underlying equations predict.
\end{abstract}

\section{Introduction}
\label{Sec:1}

In recent years, many researchers have used nonlinear partial differential equations (NPDEs) to model problems in various fields of applied science and engineering, such as chemical physics, optical fibre, plasma physics, fluid dynamics and many more\cite{1,2}. It is also significant to look at the exact explicit solutions of these models for understanding the physical phenomenon. Several methods are available to solve these types of NPDEs, such as the Lie symmetry method \cite{3,4,5,6,7}, Hirota bilinear method \cite{8,9}, B$\ddot{\text{a}}$cklund transformations method \cite{10}, inverse scattering method \cite{11} and homogeneous balance
method \cite{12}.\\\indent The Lie symmetry method is the most sophisticated to solve the NPDEs because it is based on the concrete mathematical framework and also it is the one of the powerful method to find symmetries of partial differential equations. Lie symmetry method was developed by Sophus Lie (1842–1899) in the latter half of the nineteenth century. For details of the technique, one can study the popular literature available in the Ref. \cite{3,4,5,6}\\
The extended KdV equation is \cite{13}.
 \begin{equation}\label{eq1.1}
	\frac{\partial u}{\partial t}+a\frac{\partial u}{\partial x}+\lambda u\frac{\partial u}{\partial x}-\alpha u^{2}\frac{\partial u}{\partial x}+\mu\frac{\partial^{3} u}{\partial x^{3}}+\beta\frac{\partial^{5} u}{\partial x^{5}}+\gamma_{1} u\frac{\partial^{3} u}{\partial x^{3}}+\gamma_{2}\frac{\partial u}{\partial x}\frac{\partial^{2} u}{\partial x^{2}}=0,
\end{equation}
where $a$, $\lambda$, $\alpha$, $\mu$, $\beta$, $\gamma_{1}$, $\gamma_{2}$ are arbitrary constants.
In the present work, the following nonlinear Gardner-Kawahara equation (NLGK) equation \cite{14,15} has been studied. The NLGK equation is formed by taking $\gamma_{1}=\gamma_{2}=0$.
\begin{equation}\label{eq1.2}
	\frac{\partial u}{\partial t}+a\frac{\partial u}{\partial x}+\lambda u\frac{\partial u}{\partial x}-\alpha u^{2}\frac{\partial u}{\partial x}+\mu\frac{\partial^{3} u}{\partial x^{3}}+\beta\frac{\partial^{5} u}{\partial x^{5}}=0,
\end{equation}
 NLGK equation has been investigated for symmetry reductions and conservation laws.
The above NLGK equation is a particular case of extended KdV equation. This equation widely occurs in magneto-acoustic waves in plasma physics and also in the shallow water waves with surface tension. Physically, the NLGK equation explains the solitary wave propagation in media. After considering $\alpha=\gamma_{1}=\gamma_{2}=0$ in the extended KdV equation, the Kawahara equation will be constructed.\\

In recent years, the NLGK equation attracted by many researchers. The exact travelling wave solutions for NLGK equation has been  investigated in Ref \cite{14} by implementing ($\frac{G^{\prime}}{G}$) method. Kurkina et al. \cite{15} have studied the stationary and soliton solutions of the NLGK equation.\\
 The conservation law plays an important role in the general theory and the analysis of certain systems. Conservation laws \cite{6,17,18} are essential for the understanding of the basic laws of nature and problems in mathematical physics. It is also identify the properties of the arbitrary NPDEs. Conservation laws are the most important and are broadly used in the study of the basic structures of NPDEs. In the field of physics and engineering, conservation laws allow the conclusion of a physical system, and there are some conservation laws such as conservation of mass and energy, conservation of natural resources, conservation of energy, mass, electric charge, or momentum, etc.
The best method to comprehend the physical processes of nonlinear ordinary differential equations is to find exact solutions (ODEs). The power series method \cite{7,20} and tanh method \cite{21,22} are the most effective methods to determining the exact solutions of nonlinear ODEs.

Our main focus is to study NLGK Eq. \eqref{eq1.2}. We apply the Lie symmetry method and use symmetry reduction to transform  NLGK equation into a nonlinear ODE. Thereafter, to find the exact solution of NLGK equation, we have used the power series solution method and tanh method. To finish this work, we used a multiplier technique to build the conservation rules of \eqref{eq1.2}.

\section{Lie symmetry analysis of Nonlinear Gardner-Kawahara equation}
{\label{NLGK:sec2}}
\large
In this section, the classical Lie symmetry method has been introduced for the NLGK equation. The one parameter ($\epsilon$) Lie group of transformation for the Eq.\eqref{eq1.2} with dependent and independent variable is given below
\begin{equation}
	\label{eq2.1}
	\begin{aligned}
		&x^{*}= x+\epsilon\xi^{1}(x,t,u) + O(\epsilon^{2}),\\
		&t^{*}= t +\epsilon\xi^{2}(x,t,u) +O(\epsilon^{2}),\\
		&u^{*}= u +\epsilon\eta(x,t,u)+O(\epsilon^{2}). \\
		       ....		
	\end{aligned}
\end{equation}
 Here $\epsilon$ is a group parameter and $\xi^{1},~\xi^{2},~ \eta$ are infinitesimals of the transformation for the independent and dependent variable $x, t, u$ respectively, which we have to be determined later. The associated vector fields corresponding to Eq. \eqref{eq1.2} can be written as:
 \begin{equation}\label{eq2.2}
 	V = \xi^{1}(x,t,u)\frac{\partial }{\partial x} + \xi^{2}(x,t,u)\frac{\partial }{\partial t} + \eta(x,t,u)\frac{\partial }{\partial u}. 
 \end{equation}
 If above Eq.\eqref{eq2.2} generates a symmetry Eq.\eqref{eq1.2} the infinitesimals generator $V$ must satisfy the following in variance criterion for Eq.\eqref{eq1.2}:
  \begin{equation}
  	\label{eq2.3}
 	Pr^{(5)}V(\Delta_{1})\vert\Delta_{1}=0,
 \end{equation}
 where 
 \begin{equation}\label{Eq 2.4}
 	\Delta_{1}=\frac{\partial u}{\partial t}+a\frac{\partial u}{\partial x}+\lambda u \frac{\partial u}{\partial x}-\alpha u^{2}\frac{\partial u}{\partial x}+\mu\frac{\partial^{3} u}{\partial x^{3}}+\beta\frac{\partial^{5} u}{\partial x^{5}}=0.
 \end{equation}
 $V$ is vectorfield. In general if $x$=($x_{1},x_{2},...,x_{q}$) is $q$ independent variable and $u=(u_{1},u_{2},...,u_{p})$ is $p$ dependent variable, then the related vector field $V$ will be:
\begin{equation}
	\label{eq2.5}
	V=\sum_{i=1}^{q}\xi^{i}(x,u)\partial_{x^{i}} + \sum_{j=1}^{p}\eta^{j}(x,u)\partial_{u^{j}}.
\end{equation} 
The method for finding nth order prolongation formula is given by:
 \begin{equation}
 	\label{eq2.6}
	Pr^{(n)}V = V + \sum_{k=1}^{p}\sum_{J} \eta_{k}^{J}(x,u^{(n)}\partial_{u_{J}^{k}},
\end{equation}
 where 
\begin{equation}
	\label{eq2.7}
	J=(j_{1},j_{2},..., j_{w}),~1\leq j_{w}\leq p,~0\leq J\leq n,~1\leq w \leq n 
\end{equation}
 and 
\begin{equation}
	\label{eq2.8}
	\eta_{k}^{J}(x,u^{(n)}) = D_{J} \bigg(\eta_{k} -\sum_{i=1}^{q}\xi^{i}u^{k}_{i})\bigg)+\sum_{i=1}^{q}\xi_{i}\frac{\partial }{\partial u_{J,i}^{k}},
\end{equation}
 where 
\begin{equation}
	\label{eq2.9}
	u_{i}^{k}=\frac{\partial u^{k}}{\partial u^{i}},u_{J,i}^{k}= \frac{\partial u_{J}^{k}}{\partial x^{i}}
\end{equation} 
and $D_{J}$ indicates the total derivative.
Thus the fifth order prolongation:
 \begin{equation}\label{eq2.10}
	\begin{aligned}
		& Pr^{(5)}V= V+ \eta^{x} \dfrac{\partial}{\partial u_{x}}+\eta^{t}\dfrac{\partial}{\partial u_{t}} +\eta^{xx}\dfrac{\partial}{\partial u_{xx}}+\eta^{xt}\dfrac{\partial}{\partial u_{xt}}+...
		\eta^{tt}\\&\dfrac{\partial}{\partial u_{tt}}+\eta^{xxx}\dfrac{\partial}{\partial u_{xxx}}+\eta^{xxt}\dfrac{\partial}{\partial u_{xxt}}+...\eta^{ttt}\dfrac{\partial}{\partial u_{ttt}}+\eta^{xxxx}\dfrac{\partial}{\partial u_{xxxx}}+\\& \eta^{xxxt}\dfrac{\partial}{\partial u_{xxxt}}+...\eta^{tttt}\dfrac{\partial}{\partial u_{tttt}}+\eta^{xxxxx}\dfrac{\partial}{\partial u_{xxxxx}}+\eta^{xxxxt}\dfrac{\partial}{\partial u_{xxxxt}}+...\\ &\eta^{ttttt}\dfrac{\partial}{\partial u_{ttttt}}.
	\end{aligned}
\end{equation}
The obtained result after applying the fifth order prolongation in Eq.\eqref{eq1.2}
\begin{equation}\label{eq2.11}
	\eta_{t}+a\eta_{x}+\lambda\eta\eta_{x}-a\eta^{2}\eta^{x}+\mu\eta^{xxx}+\beta\eta^{xxxxx}=0.
\end{equation} 
One can solve the Eq.\eqref{eq2.11} by substituting the coefficient of several monomials equal to zero, leading to a system of determining equations that have been solved with MAPLE software \cite{shingareva2011,kumar2014}.
\begin{equation}
	\label{eq2.12}
	\xi^{1}=c_{2},\xi^{2}=c_{1},\eta=0
\end{equation} 
The infinitesimals of Eq.\eqref{eq1.2} is given by:
\begin{equation}
	\label{eq2.13}
	X_{1}= \frac{\partial }{\partial t}, X_{2}=\frac{\partial}{\partial x}
\end{equation}
The infinitesimals generator of Eq.\eqref{eq1.2} will be:
\begin{equation}\label{eq2.14}
	V = c_{1}X_{1}+c_{2}X_{2}
\end{equation} 
The two translation symmetries of Eq.\eqref{eq2.13} used to find the symmetry reductions of Eq.\eqref{eq1.2}. Here we find a linear combination of these Lie symmetries, $X$= $X_{1}+ cX_{2}$, where $c$ is constant, has been used. The Lagrange equation will be:
 \begin{equation}
 	\label{eq2.15}
	\dfrac{dx}{1}=\dfrac{dt}{c}=\dfrac{du}{0}
\end{equation}
The invariant solution will be
\begin{equation}\label{eq2.16}
	\begin{aligned}
		& z=(x-ct)\\
		&  U=u
	\end{aligned}
\end{equation}
 Now the group-invariant solution becomes $u(x,t)$ = $U(z)$. Now $U$ and $z$ will be the new dependent and independent variables. After applying invariant solution into Eq.\eqref{eq1.2} reduces the following fifth-order non-linear ODE:
\begin{equation}\label{eq2.17}
	\beta U^{\prime\prime\prime\prime\prime}+\mu U^{\prime\prime\prime}-\alpha U^2U^{\prime}+\lambda U U^{\prime}+aU^{\prime}- c U^{\prime}=0.
\end{equation}
 \section{Exact explicit solution by power series method}
{\label{NLGK:sec3}}
\large
In general, it is not an easy task to find exact and analytic solutions for this type of nonlinear ODE Eq.\eqref{eq2.17}, by applying the elementary functions and integrals. In this section, we will derive the exact solution of the reduced Eq.\eqref{eq2.17} by using the power series method. The power series method \cite{7,20} is a most significant method to solve higher-order ODEs. \\
Consider, 
\begin{equation}\label{Eq 20}
	U(z)= \sum_{n=0}^{\infty}c_{n}z^{n}.
\end{equation}
Differentiating Eq.\eqref{Eq 20} and substituting into Eq.\eqref{eq2.17} leads to

\begin{eqnarray}\label{Eq 21}
	c\sum_{n=0}^{\infty} (n+1)c_{n+1}z^{n}-a\sum_{n=0}^{\infty}(n+1)c_{n+1}z^{n}-\lambda\sum_{n=0}^{\infty} c_{n}z^n\sum_{n=0}^{\infty}(n+1)c_{n+1}z^{n}+\\\nonumber
	\alpha \sum_{n=0}^{\infty} c_{n}z^{n}\sum_{n=0}^{\infty}c_{n}z^{n}\sum_{n=0}^{\infty}(n+1)c_{n+1}z^{n}- \mu\sum_{n=0}^{\infty}(n+1)(n+2)(n+3)c_{n+3}z^{n}-\\
	\beta\sum_{n=0}^{\infty}(n+1)(n+2)(n+3)(n+4)(n+5)c_{n+5}z^{n}=0,\nonumber
\end{eqnarray}
which implies 
\begin{eqnarray}\label{22}
	c\sum_{n=1}^{\infty}(n+1)c_{n+1}z^{n}+cc_{0}-a\sum_{n=1}^{\infty}(n+1)c_{n+1}z^{n}-ac_{1}-\lambda\sum_{n=1}^{\infty}c_{n}z^{n}\sum_{n=1}^{\infty}(n+1)c_{n+1}z^{n}-\lambda c_{0}c_{1}+\\\nonumber
	\alpha\sum_{n=1}^{\infty}c_{n}z^{n}\sum_{n=1}^{\infty} c_{n}z^{n}\sum_{n=1}^{\infty}(n+1)c_{n+1}z^{n}-\alpha c_{0}^{2}c_{1}-\mu\sum_{n=1}^{\infty}(n+1)(n+2)(n+3)c_{n+3}z^{n}-6\mu c_{3}-\\
	\beta\sum_{n=1}^{\infty}(n+1)(n+2)(n+3)(n+4)(n+5)c_{n+5}z^{n}-120\beta c_{5}=0,\nonumber
\end{eqnarray}
On elaborating the above Eq.
\begin{eqnarray}\label{23}
	c\sum_{n=1}^{\infty}(n+1)c_{n+1}z^{n}-a\sum_{n=1}^{\infty}(n+1)c_{n+1}z^{n}-\lambda\sum_{n=1}^{\infty}\sum_{k=0}^{n}(n-k+1)c_{k}c_{n-k+1}z^{n}\\\nonumber
	+\alpha\sum_{n=1}^{\infty}\sum_{k=0}^{n}\sum_{i=0}^{k}(n-k+1)c_{n-k+1}c_{i}c_{k-i}z^{n}-\mu\sum_{n=1}^{\infty}(n+1)(n+2)(n+3)c_{n+3}z^{n}-\\\nonumber
	\beta\sum_{n=1}^{\infty}(n+1)(n+2)(n+3)(n+4)(n+5)c_{n+5}z^{n}+\\\nonumber
	{cc_{0}-ac_{1}-\lambda c_{0}c_{1}-6\mu c_{3}-120\beta c_{5}-\alpha c_{0}^{2} c_{1}}=0
\end{eqnarray}
Collecting the same powers of $ z $, one can obtain:
\begin{eqnarray}\label{Eq 24}
	\sum_{n=1}^{\infty}[c(n+1)c_{n+1}-a(n+1)c_{n+1}-\sum_{k=0}^{n}\lambda(n-k+1)c_{k}c_{n-k+1}\\\nonumber
	+\sum_{k=0}^{n}\sum_{i=0}^{k}(n-k+1)c_{n-k+1}c_{i}c_{k-i}-\mu(n+1)(n+2)(n+3)c_{n+3}\\\nonumber- 
	\beta(n+1)(n+2)(n+3)(n+4)(n+5)c_{n+5}]z^{n}\\+{cc_{0}-ac_{1}-\lambda c_{0}c_{1}-6\mu c_{3}-120\beta c_{5}-\alpha c_{0}^{2}c_{1}}=0.
	\nonumber
\end{eqnarray}
On solving Eq.\eqref{Eq 24},
\begin{equation}\label{25}
	c_{5}=\dfrac{cc_{0}-ac_{1}-\lambda c_{0}c_{1}-6\mu c_{3}-\alpha c_{0}^2 c_{1}}{120 \beta},
\end{equation}
\begin{multline}\label{26}
	\begin{aligned}
		c_{n+5}=&\dfrac{1}{(n+1)(n+2)(n+3)(n+4)(n+5) \beta}\bigg[c(n+1)c_{n+1}-a(n+1)c_{n+1}\\&-\sum_{k=0}^{n}\lambda (n-k+1)c_{k}c_{n-k+1}+\sum_{k=0}^{n}\sum_{i=0}^{k}(n-k+1)c_{n-k+1}c_{i}c_{k-i}\\&-\mu (n+1)(n+2)(n+3)c_{n+3}\bigg],
	\end{aligned}
\end{multline} 
On the other hand, for finding the coefficient $ c_{n}(n\geq5)$ of the power series Eq.\eqref{25} and Eq.\eqref{26} has been determined by taking arbitrary constant $c,c_{0}, c_{1}, c_{3}, c_{2}, c_{4}, a, \alpha, \lambda,  \beta, \mu$
and also clearly, the power series solution depends on  Eq.\eqref{25} and Eq.\eqref{26}.
 \begin{equation}
	U(z)=c_{0}+c_{1}z+c_{2}z^{2}+c_{3}z^{3}+c_{4}z^{4}+c_{5}z^{5}+\sum_{n=1}^{\infty}c_{n+5}z^{n+5} .
\end{equation}
Furthermore, we obtain the explicit power series solution for Eq.\eqref{eq1.2}
\begin{multline}
	\begin{aligned}
		u(x,t)= & c_{0}+c_{1}(x-ct)+c_{2}(x-ct)^{2} +c_{3}(x-ct)^{3} +c_{4}(x-ct)^{4} \\&+ \bigg( \dfrac{cc_{0}-ac_{1}-\lambda c_{0}c_{1}-6\mu c_{3}-\alpha c_{0}^2 c_{1}}{120 \beta} \bigg)(x-ct)^{5}\\&+
		\sum_{n=1}^{\infty} \bigg( \dfrac{1}{(n+1)(n+2)(n+3)(n+4)(n+5) \beta}\bigg[c(n+1)c_{n+1}\\&-a(n+1)c_{n+1}-\sum_{k=0}^{n}\lambda (n-k+1)c_{k}c_{n-k+1}\\&+\sum_{k=0}^{n}\sum_{i=0}^{k}(n-k+1)c_{n-k+1}c_{i}c_{k-i}\\&-\mu (n+1)(n+2)(n+3)c_{n+3}\bigg] \bigg)(x-ct)^{n+5}.
	\end{aligned}
\end{multline}
\section{Exact explicit solution by tanh method}
This section consists the algorithm of tanh method \cite{21,22} and its implementation for solving NLGK equation.
\subsection{\textbf{Algorithm for tanh method}}
\large
The main steps of tanh method  are discussed as follows\\
\textbf{Step 1}: Let us assume a NPDEs in this form
\begin{equation}\label{29}
	F(u,u_{x},u_{xx},...,u_{t})=0
\end{equation} where, $u=u(x,t)$ is an unknown function. \\
\textbf{Step 2}: Using the transformation,
\begin{equation}\label{30}
	u(x,t)=U(z),~~ z=(x-ct),
\end{equation} where, $ c $ is constant to be determined later, the PDEs \eqref{29} is reduced to the ODEs as follows
\begin{equation}\label{31}
	F(U,U^{\prime},U^{\prime\prime},...,-cU^{\prime})=0  
\end{equation}
\textbf{Step 3}: Suppose the solution of Eq.\eqref{31} can be expressed in the form
\begin{equation}\label{32}
	U(z)= a_{0}+\sum_{i=1}^{M}a_{i}Y^{i}
\end{equation} where $ a_{i} $ are constants to be determined, the integer $M$ can be calculated by balancing the highest order derivative term with the highest order nonlinear term occuring in Eq.\eqref{31} and $Y=\text{tanh}(z)$ is new independent variable. The derivatives of Eq.\eqref{32} can be obtained as follows
 \begin{equation}\label{33}
	\begin{aligned}  
		& \frac{dU}{dz}\longrightarrow (1-Y^2)\frac{dU}{dY}\\
		&  \frac{d^2 U}{dz^2}\longrightarrow (1-Y^2)\bigg(-2Y\frac{dU}{dY}+(1-Y^2)\frac{d^2U}{dY^2}\bigg) 
	\end{aligned}
\end{equation} and so on for higher.\\
\textbf{Step 4} Substituting Eq.\eqref{32} and \eqref{33} in Eq.\eqref{31}, the solution will get in terms of $Y^{i}(i=0,1,2,3...).$ Then equating to zero of same degree $Y^{i}(i=0,1,2,3...)$ results into a set of algebraic equations for $a_{i}(i=0,1,2...),c, a, \alpha,\beta,\mu .$\\
\textbf{Step 5}: Solving the resulting algebraic system, the exact explicit solution of Eq.\eqref{29} can be obtained.
 \subsection{\textbf{Implementation of tanh method for solving NLGK equation }}
\large
For finding the exact solution of Eq.\eqref{eq2.17}, the above algorithm is implemented.\\
Now integrating the Eq.\eqref{eq2.17} with respect to z, we will get
\begin{equation}\label{34}
	(a-c) U+\dfrac{\lambda}{2} U^2-\dfrac{\alpha}{3} U^3+\mu U^{\prime\prime\prime}+\beta U^{\prime\prime\prime\prime\prime}=0.  
\end{equation}
Hence M=2. From Eq.\eqref{32}, the solution Eq.\eqref{eq2.17} can be expressed as 
\begin{equation}\label{35}
	U(z)=a_{0}+a_{1}Y+a_{2}Y^2
\end{equation}
where, $ Y=\text{tanh}(z) $. By substituting Eq.\eqref{35} and \eqref{33} into the Eq.\eqref{34}, collecting all terms of same powers of $ Y $ and equating them to zero, a set of algebraic equations can be obtained. On solving the resulting algebraic system, the following set of solutions can be achieved.\\
 \textbf{Set 1}
\begin{equation}\label{36}
	\begin{aligned}
		&a_{0}=0,~a_{1}=0,~c=\dfrac{5a+\alpha a_{2}^2}{5},~\beta=\dfrac{\alpha a_{2}^2}{360},~\lambda=\dfrac{16 \alpha a_{2}}{15},~\mu=\dfrac{\alpha a_{2}^2}{45}\\
		&u(x,t)=a_{2} \bigg(\text{tanh}\bigg[x-\bigg(\dfrac{5a+\alpha a_{2}^2}{5}\bigg)t\bigg]\bigg)^2 
	\end{aligned}
\end{equation}\\
\textbf{Set 2}
\begin{equation}\label{37}
	\begin{aligned}
		& a_{1}=0,~a_{2}=-2a_{0},~\alpha=\dfrac{15(a-c)}{7 a_{0}^2},~\beta=\dfrac{a-c}{42},~\lambda=\dfrac{4(a-c)}{7 a_{0}},~\mu=\dfrac{(a-c)}{3}\\
		& u(x,t)= a_{0}-2a_{0}\bigg(\text{tanh}[x-ct]\bigg)^2
	\end{aligned}
\end{equation}\\
\textbf{Set 3}
\begin{equation}\label{38}
	\begin{aligned}
		&a_{1}=0,~a_{2}=-a_{0},~c=\dfrac{(45a+15\lambda a_{0}-8\alpha a_{0}^2)}{45},\\&~\beta=\dfrac{\alpha a_{0}^2}{360},~
		\mu=\dfrac{3\lambda a_{0}-2\alpha a_{0}^2}{36} \\
		&u(x,t)= a_{0}-a_{0}\bigg(\text{tanh}\bigg[x-\bigg(\dfrac{45a+15\lambda a_{0}-8\alpha a_{0}^2}{45}\bigg)t\bigg]\bigg)^2
	\end{aligned}
\end{equation}
\begin{figure}[th!]
	\begin{minipage}[b]{0.4\textwidth}
		\includegraphics[width=\textwidth]{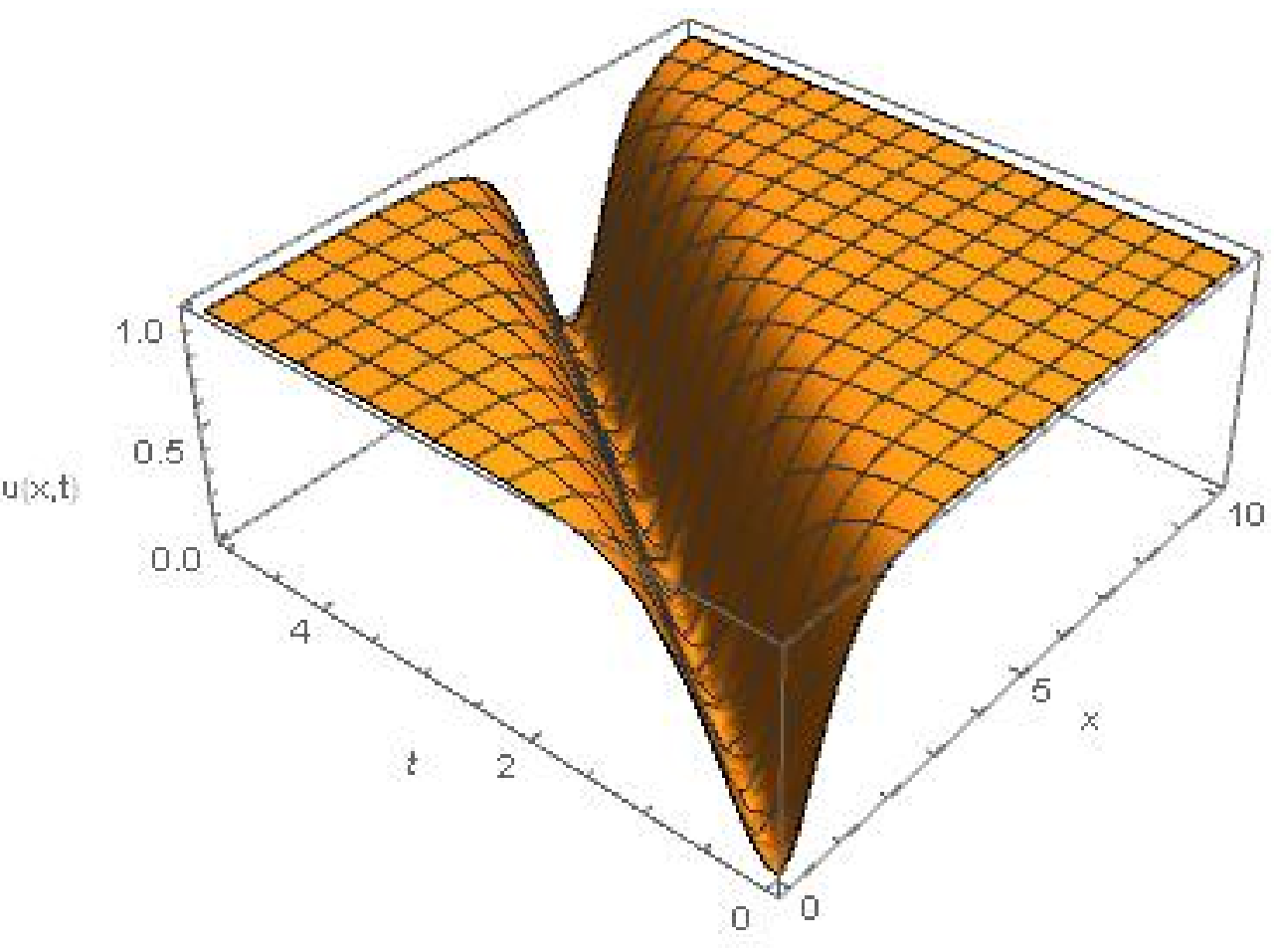}
		\caption{Soliton solution of Eq. \eqref{36}, by setting parameters $a_{2}=1,~a=1,~\alpha=2$.}
	\end{minipage}
	\hfill
	\begin{minipage}[b]{0.4\textwidth}
		\includegraphics[width=\textwidth]{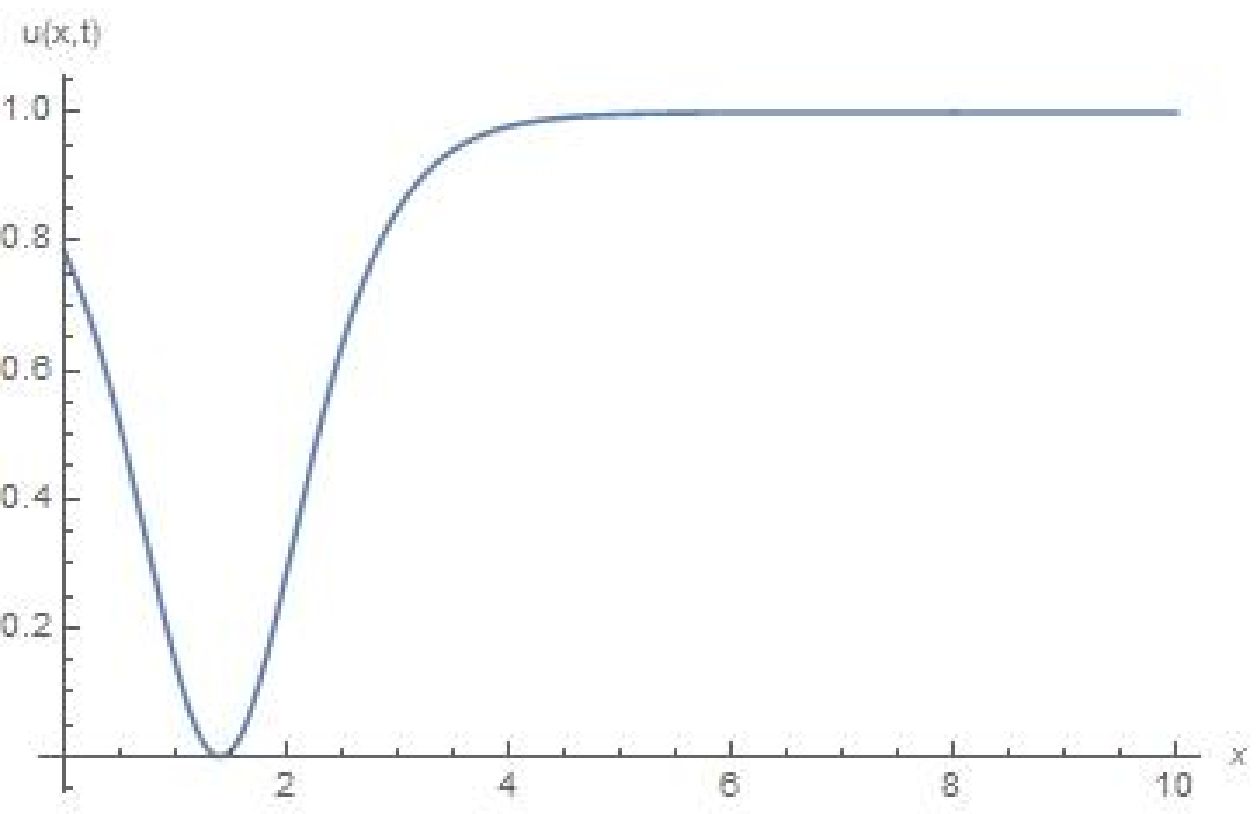}
		\caption{2D plot Eq.\eqref{36}, by setting suitable arbitrary parameters $a_{2}=1,~a=1,~\alpha=2,~t=1$.}
	\end{minipage}
\end{figure}

\begin{figure}[th!]
	\begin{minipage}[t]{0.4\textwidth}
		\includegraphics[width=\textwidth]{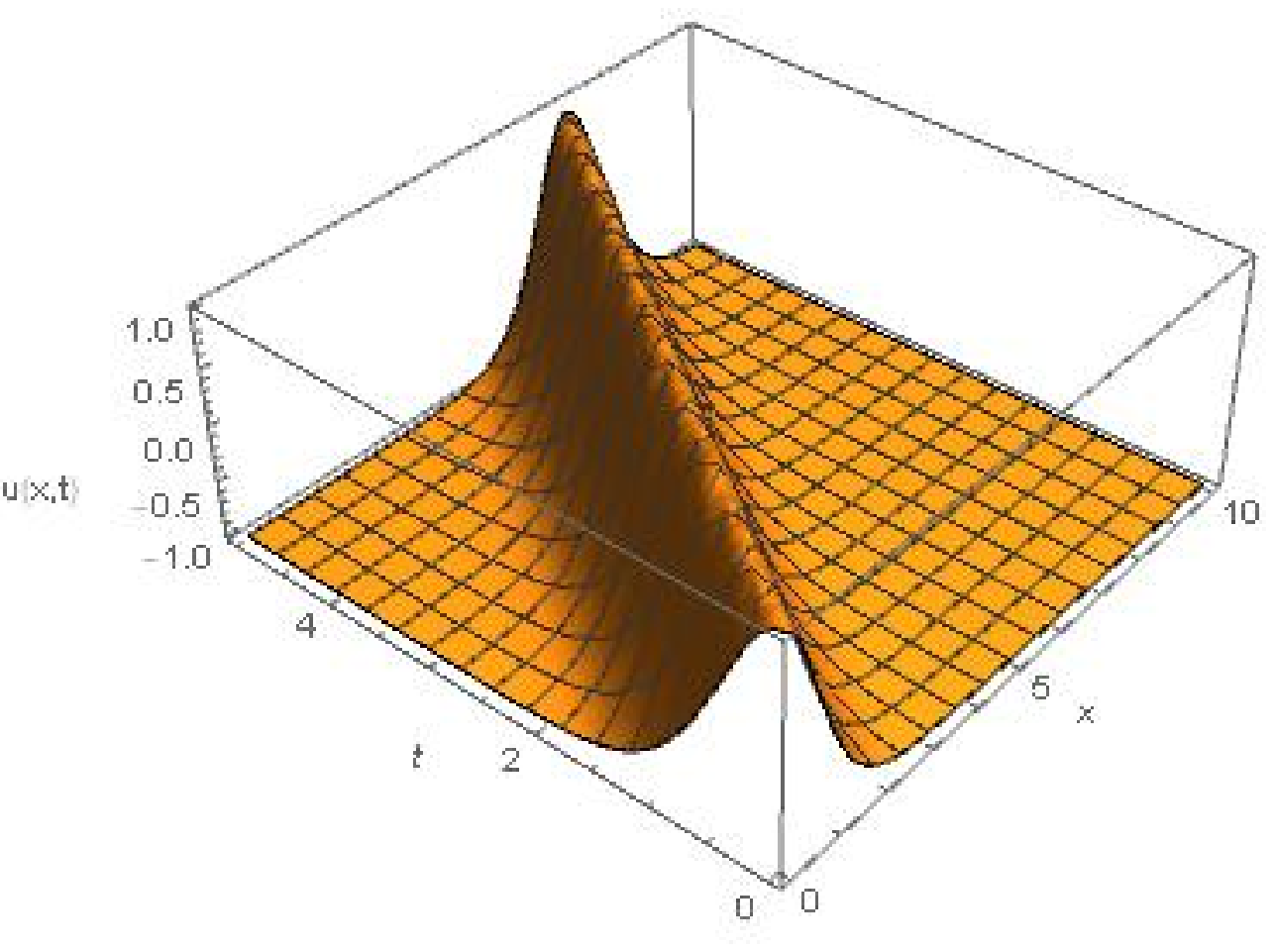}
		\caption{Solitary wave profile of Eq.\eqref{37} by setting parameters $a_{0}=1,~c=1.4$. }
	\end{minipage}
	\hfill
	\begin{minipage}[b]{0.4\textwidth}
		\includegraphics[width=\textwidth]{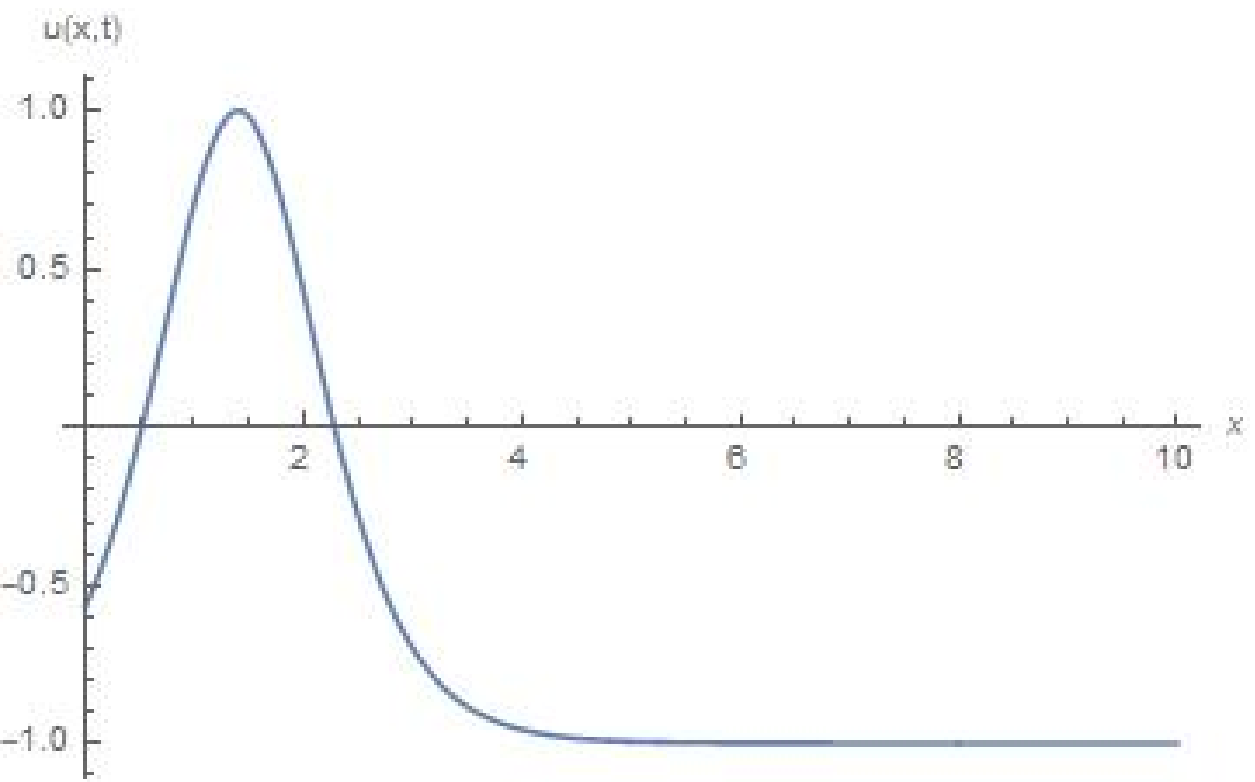}
		\caption{2D plot of Eq.\eqref{37} by setting parameters $a_{0}=1,~c=1.4,~t=1$.}
	\end{minipage}
\end{figure}

\begin{figure}[th!]
	\begin{minipage}[b]{0.4\textwidth}
		\includegraphics[width=\textwidth]{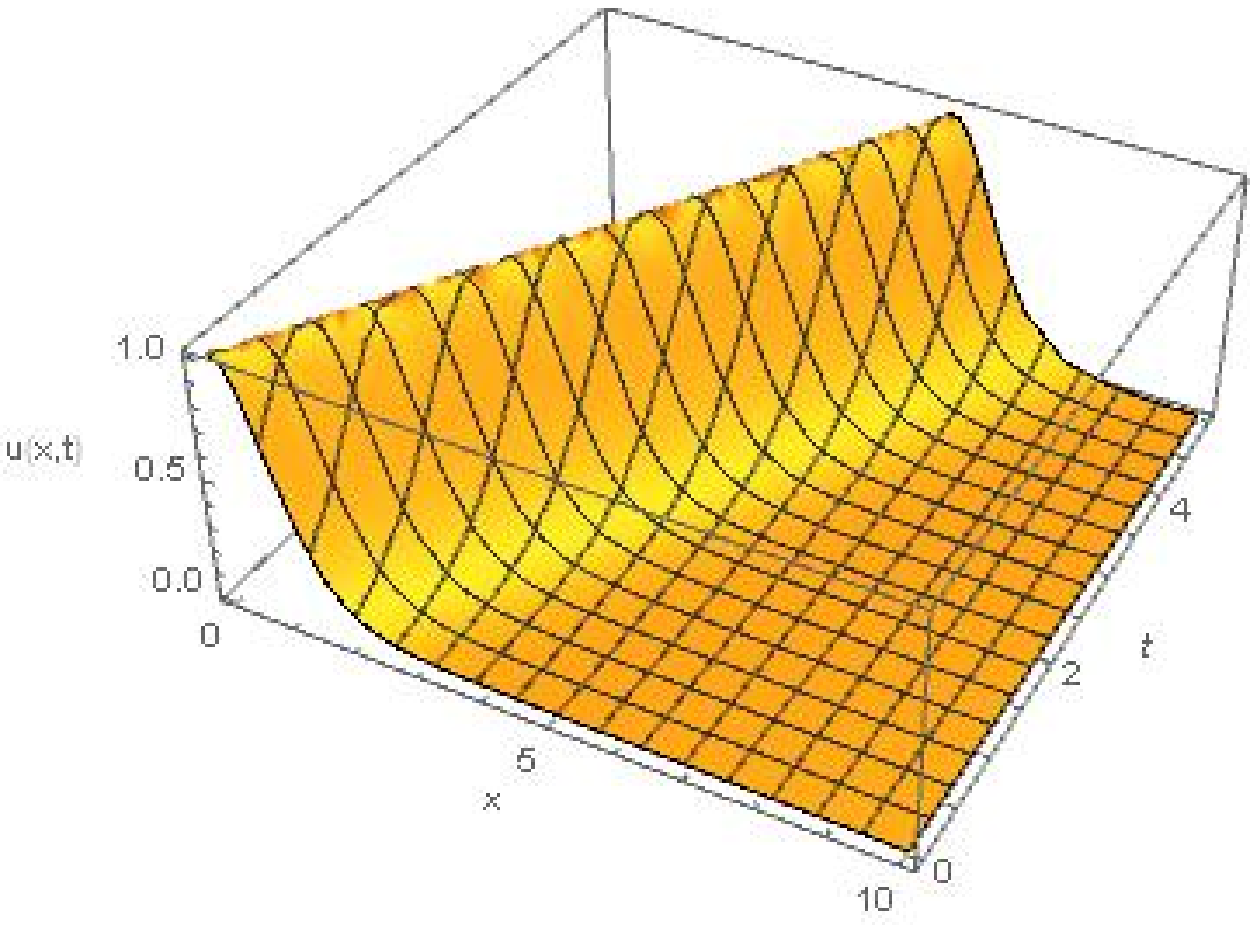}
		\caption{Solitary wave profile of Eq.\eqref{38} by setting    all parameters $a_{0}, a,\lambda$ to unity.}
	\end{minipage}
	\hfill
	\begin{minipage}[b]{0.4\textwidth}
		\includegraphics[width=\textwidth]{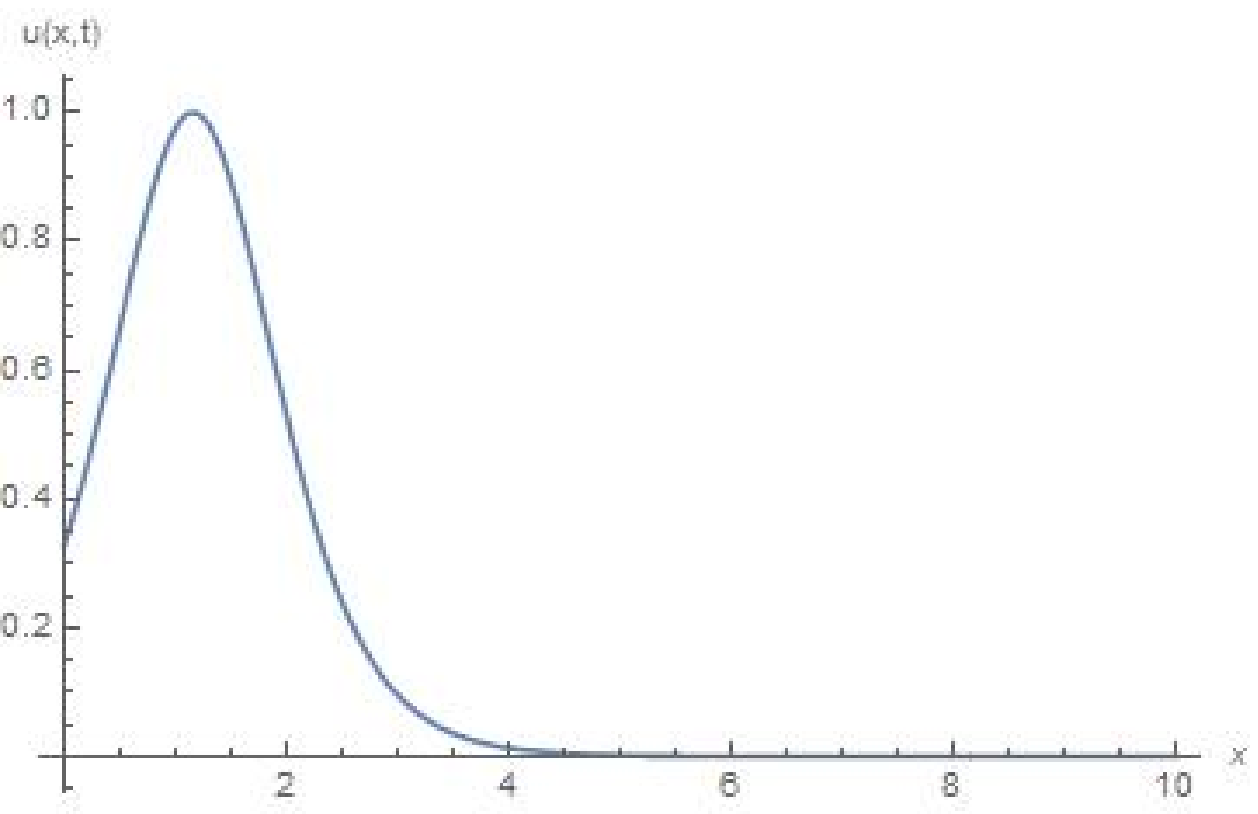}
		\caption{2D plot of Eq.\eqref{38} by setting  all parameters $a_{0}, a,\lambda,t$ to unity.}
	\end{minipage}
\end{figure}
\section{Construction of conservation laws via multiplier method} 
{\label{NLGK:sec4}}
\large
In this section, non-trival local conservation laws for the NLGK equation have been constructed. The multiplier method of Anco and Bluman \cite{6,17,18} has been used for this construction. The detailed application of this method can be seen in recent work \cite{23,24,25}. \\
{\textbf{Definition 1}}
The Euler operator for the dependent variable $u^{j}$ is  defined by
\begin{equation}\label{Eq 39}
	E_{u^{j}}=\dfrac{\delta}{\delta u^{j}}=\dfrac{\partial}{\partial u^{j}} + \sum_{p=1}^{\infty}(-1)^{p}{D_{i_{1}}}...{D_{i_{p}}} \dfrac{\partial}{\partial u_{{i_{1}...i_{p}}}^{j}}
\end{equation} for each $j=1,...,m.$\\
{\textbf{Definition 2}}
The total differentiation which is defined as,
 \begin{equation}\label{Eq 40}
	D_{i}= \dfrac{\partial}{\partial x^{i}}+ u_{i}\dfrac{\partial}{\partial u}+ u_{ij}\dfrac{\partial}{\partial u_{j}}+...
\end{equation}
which is defined with respect to independent variable $x$=($x^{1}$ ,...,$x^{n}$). The Euler operator Eq.\eqref{Eq 39} which can annihilate any divergence expression $D_{i}$ $\Psi^{i}$(u)        \cite{singh2018}.
This method investigates determination of the zero order multiplier $\Lambda(x,t,u)$ such that \begin{equation}\label{Eq 41}
	E_{u}[\Lambda(x,t,u)(u_{t}+au_{x}+\lambda uu_{x}-\alpha u^{2}u_{x}+
	\mu u_{xxx}+\beta u_{xxxxx}]=0
\end{equation}
where, $E_{u}$ is the Euler operator defined in Eq.\eqref{Eq 39}. On expanding the Eq.\eqref{Eq 41}, the following determination system can be obtained for multiplier $\Lambda(x,t,u)$
\begin{multline}\label{Eq 42}
	\Lambda_{uu}=0,~\Lambda_{xu}=0,~
	-\Lambda_{t}-\beta\Lambda_{xxxxx}-\mu \Lambda_{xxx}-\Lambda_{x}(-\alpha u^{2} +\lambda u +a)=0
\end{multline} By solving above equation \eqref{Eq 42} one can obtain
\begin{equation}\label{Eq 43}
	\Lambda_{1}(x,t,u)=1, ~\Lambda_{2}(x,t,u)=u.
\end{equation} Each multiplier from Eq.\eqref{Eq 43} determines local conservation laws in the format \begin{equation}\label{44}
	D_{x}\Psi_{1}(x,t,u)+D_{t}\Psi_{2}(x,t,u)=0
\end{equation} with the characteristics form :
\begin{equation}\label{Eq 45}
	\begin{aligned}
		& \Lambda(x,t,u)(u_{t}+au_{x}+\lambda uu_{x}-\alpha u^{2}u_{x}+
		\mu u_{xxx}+\beta u_{xxxxx})\\&= D_{x}\Psi_{1}(x,t,u)+D_{t}\Psi_{2}(x,t,u)
	\end{aligned}
\end{equation} The inversion of divergence expression Eq.\eqref{Eq 45} can be carried out by 2-dimensional homotopy operator \cite{singh2018} and results read the following:
\begin{multline}
	\begin{aligned}
		\bullet \Lambda_{1}(x,t,u)=1\\
		&\psi_{1}(x,t,u)=\dfrac{-1}{3}\alpha u^3+\dfrac{1}{2}\lambda u^2+au+\beta u_{xxxx}+\mu u_{xx},\\
		& \psi_{2}(x,t,u)= u\\
		\bullet \Lambda_{2}(x,t,u)=u\\
		&\psi_{1}(x,t,u)=\dfrac{-1}{4}\alpha u^{4}+\dfrac{1}{3}\lambda u^{3} +\dfrac{1}{2}a u^{2}+\beta u u_{xxxx}-\beta u_{x} u_{xxx}+\\&\dfrac{1}{2} \beta u_{xx}^{2}+\mu u u_{x}-\dfrac{1}{2} \mu u_{x}^{2},\\
		&\psi_{2}(x,t,u)= \dfrac{u^2}{2}
	\end{aligned}
\end{multline}
The non-triviality of conservation laws can be checked from the fact that divergence expression $D_{x}\Psi_{1}(x,t,u)$+$D_{x}\Psi_{2}(x,t,u)$ vanishes on solution space of the NLGK equation, and most importantly, none of them vanishes on solution space of NLGK equation.
 \section{Conclusions}\label{sec5}
\large
In the present work, Lie symmetry analysis on the nonlinear NLGK equation has been studied. The Lie group method is utilized to acquire the symmetry reductions of NLGK equation. Using the Lie symmetry reductions and the power series method, the exact analytic solution of NLGK equation has been derived from the reduction equation. Also, the exact solutions are obtained by implementing the tanh method to the reduction equation. From the analysis of the acquired solutions, it may be concluded that NLGK equation produces a collection of travelling wave solutions. Finally, the conservation laws for NLGK equation have been also constructed by utilizing the multiplier method.
\section*{Acknowledgment}
The authors are grateful to the referees for their educative and constructive feedback, which helped to improve the manuscript. The first author would like to thank the Department of Science and Technology (DST) of Goverment of Odisha, India for financial support in the form of PhD sponsorship.

\end{document}